\documentclass[12pt]{amsart}

\usepackage{amsfonts,amsmath,latexsym,amssymb,verbatim,amsbsy}
\usepackage{amsthm, graphicx, times}

\def\vu{U}

\def\vk{k}

\def\cF{{\mathcal F}}
\def\vf{\varphi}

\def\cH{{\mathcal H}}
\def\D{\Delta}
\def\L{Q}
\def\wt{\widetilde}

\def\Supp{{\rm Supp }}
\def\eqdefa{\buildrel\hbox{\footnotesize def}\over =}

\def \vanhel {B^{2/3}_{3, c(\mathbb N)}}

\newcommand{\RR}{{\mathbb R}}

\newcommand{\ZZ}{{\mathbb Z}}
\newcommand{\NN}{{\mathbb N}}
\newcommand{\CC}{{\mathbb C}}

\newcommand{\Sum}{\displaystyle \sum}

\newcommand{\be}{\begin{equation}}
\newcommand{\ee}{\end{equation}}
\newcommand{\la}{\label}
\newcommand{\ba}{\begin{array}{c}}
\newcommand{\ea}{\end{array}}

\def \o {\omega}
\def \O {\Omega}

\newcommand{\tr}{{\rm Tr}}

\def\p{\partial}

\def \< {\langle}
\def \> {\rangle}
\def \p {\partial}
\def \ra {\rightarrow}

\def \van {B^{1/3}_{3, c(\mathbb N)}}
\def \bes {B^{1/3}_{3,\infty}}
\def \l {\lambda}

\def \n {\nabla}

\newtheorem{thm}{Theorem}[section]
\newtheorem{prop}[thm]{Proposition}
\newtheorem{lemma}[thm]{Lemma}

\newtheorem{defi}[thm]{Definition}
\newtheorem{cor}[thm]{Corollary}

\begin{document}

\title[Energy conservation]{Energy conservation and Onsager's conjecture for the Euler equations}
\author{ A. Cheskidov}
\address[A. Cheskidov]
         {Department of Mathematics\\
          University of Michigan\\
          Ann Arbor, MI 48109}
\email{acheskid@umich.edu}

\author{P. Constantin}
\address[P. Constantin]
         {Department of Mathematics\\
          University of Chicago\\
           Chicago, IL 60637}
\email{const@cs.uchicago.edu}

\author{S. Friedlander}
\address[S. Friedlander and R. Shvydkoy]
{Department of Mathematics, Stat. and Comp. Sci.\\
        University of Illinois\\
        Chicago, IL 60607}
\email{susan@math.northwestern.edu}

\author{R. Shvydkoy}
\email{shvydkoy@math.uic.edu}

\date{April 5, 2007}

\begin{abstract}
Onsager conjectured that weak solutions of the Euler equations for
incompressible fluids in $\RR^3$ conserve energy only if they have a
certain minimal smoothness, (of order of $1/3$ fractional
derivatives) and that they dissipate energy if they are rougher. In
this paper we prove that energy is conserved for velocities in the
function space $B^{1/3}_{3,c(\NN)}$. We show that this space is
sharp in a natural sense. We phrase the energy spectrum in terms of
the Littlewood-Paley decomposition and show that the energy flux is
controlled by local interactions. This locality is
shown to hold also for the helicity flux; moreover, every weak
solution of the Euler equations that belongs to $B^{2/3}_{3,c(\NN)}$
conserves helicity. In contrast, in two dimensions, the strong
locality of the enstrophy holds only in the ultraviolet
range.
\end{abstract}



\maketitle

\section{Introduction}
The Euler equations for the motion of an incompressible inviscid fluid are
\begin{equation} \label{e:Euler}
\frac{\partial u}{\partial t} + (u \cdot \nabla)u = - \nabla p,
\end{equation}
\begin{equation}
\nabla \cdot u =0,
\end{equation}
where $u(x,t)$ denotes the $d$-dimensional velocity, $p(x,t)$ denotes the pressure, and $x\in\RR^d$. We mainly consider the case $d=3$.
When $u(x,t)$ is a classical solution, it follows directly that
 the total energy $E(t) = \frac{1}{2} \int |u|^2 \, dx$ is conserved. However, conservation of energy may fail for weak solutions (see Scheffer \cite{Shn},
Shnirelman \cite{Sch}). This possibility has given rise to a considerable body of literature and it is closely connected with statistical theories of turbulence envisioned 60 years ago by Kolmogorov and Onsager. For reviews see, for example, Eyink and Sreenivasan \cite{ES}, Robert \cite{R}, and Frisch \cite{F}.

Onsager \cite{O} conjectured that in 3-dimensional turbulent flows, energy dissipation might exist even in the limit of vanishing viscosity. He suggested that an appropriate mathematical description of turbulent flows (in the inviscid limit) might be given by weak solutions  of the Euler equations that are not regular enough to conserve energy. According to this view, non-conservation of energy in a turbulent flow might occur not only from viscous dissipation, but also from lack of smoothness of the velocity. Specifically, Onsager conjectured that weak solutions of the Euler equation with H\"{o}lder  continuity exponent $h>1/3$ do conserve energy and that turbulent or anomalous dissipation occurs when $h\leq 1/3$. Eyink \cite{E} proved energy conservation under a stronger assumption. Subsequently, Constantin, E and Titi \cite{cet} proved energy conservation for $u$ in the Besov space
$B^{\alpha}_{3,\infty}$, $\alpha>1/3$. More recently the result was proved under a
slightly weaker assumption by Duchon and Robert \cite{DR}.

In this paper we sharpen the result of \cite{cet}: we prove that energy is conserved for velocities in the Besov space of tempered distributions $B^{1/3}_{3,p}$. In
fact we prove the result for velocities in the slightly larger space $B^{1/3}_{3,c(\NN)}$ (see Section 3). This is a space in which the ``H\"{o}lder exponent'' is exactly $1/3$, but the slightly better regularity is encoded in the summability condition.
The method of proof combines the approach of \cite{cet} in bounding the trilinear term in (3) with a suitable choice of the test function for weak solutions in terms of a Littlewood-Paley decomposition. Certain cancelations in the trilinear term become apparent using this decomposition.
We observe that the space $B^{1/3}_{3,c(\NN)}$ is sharp in the context of no anomalous dissipation. We give an example of a divergence free vector field in $B^{1/3}_{3,\infty}$ for which the energy flux due to the trilinear term is bounded
from below by a positive constant. This construction follows ideas in \cite{E}. However, because it is not a solution of the unforced
Euler equation, the example does not prove that indeed there exist unforced solutions to the Euler equation that live in
$B^{1/3}_{3,\infty}$ and dissipate energy.

Experiments and numerical simulations indicate that for many turbulent flows the energy dissipation rate appears to remain positive at large Reynolds numbers. However, there are no known
 rigorous lower bounds for  slightly viscous Navier-Stokes equations. The existence of a weak solution of Euler's equation,
with positive smoothness and that does not conserve energy remains an open question. For a discussion see, for example, Duchon and Robert \cite{DR}, Eyink \cite{E}, Shnirelman \cite{Shn},
Scheffer \cite{Sch}, de Lellis and Szekelyhidi \cite{LS}.

We  note that the proof in Section 3 applied to Burger's equation for 1-dimensional compressible flow gives conservation of energy in $B^{1/3}_{3,c(\NN)}$. In this case it is easy to show that conservation of energy can fail in $B^{1/3}_{3,\infty}$ which is the sharp space for shocks.

The Littlewood-Paley approach to the
issue of energy conservation versus turbulent dissipation is mirrored in a
study of a discrete dyadic model for the forced Euler equations \cite{CFP1, CFP2}.
By construction, all the interactions in that  model system are local
and energy cascades strictly to higher wave numbers. There is a
unique fixed point which is an exponential global attractor.
Onsager's conjecture is confirmed for the model in both directions,
i.e. solutions with bounded $H^{5/6}$ norm satisfy the energy balance
condition and turbulent dissipation occurs  for all solutions when
the $H^{5/6}$ norm becomes unbounded, which happens in finite time. The absence of anomalous dissipation for inviscid
shell models has been obtained in \cite{clt} in a space with regularity logarithmically higher than $1/3$.

In Section~\ref{s:EnergyFlux} we present the definition of the energy flux employed in the
paper. This is the flux of the Littlewood-Paley spectrum, (\cite{C97}) which is
a mathematically convenient variant of the physical concept of flux from the
turbulence literature. Our estimates employing the Littlewood-Paley decomposition produce not only a sharpening of the conditions under which there is no anomalous dissipation, but also provide detailed information concerning the cascade of energy flux through frequency space. In section 3.3. we prove that the energy flux through the sphere of radius $\kappa$ is controlled primarily by scales of order $\kappa$. Thus we give a mathematical justification for the physical intuition underlying much of turbulence theory, namely that the flux is controlled by local interactions (see, for example, Kolmogorov \cite{K41} and also \cite{E05}, where sufficient conditions for locality were described). Our analysis makes precise an exponential decay of nonlocal contributions to the flux that was conjectured by Kraichnan \cite{Kraichnan59}.

The energy is not the only scalar quantity that is conserved under evolution by classical solutions of the Euler equations. For 3-dimensional flows the helicity is an important quantity related to the topological configurations of vortex tubes (see, for example, Moffatt and Tsinober \cite{MT}). The total helicity is conserved for smooth ideal flows. In Section 4 we observe that the techniques used in Section 3 carry over exactly to considerations of the helicity flux, i.e.,  there is locality for turbulent cascades of helicity and every weak solution of the Euler equation that belongs to $B^{2/3}_{3,c(\NN)}$ conserves helicity. This strengthens a recent result of Chae \cite{Chae}. Once again our argument is sharp in the sense that a divergence free vector field in $B^{2/3}_{3,\infty}$ can be constructed to produce an example for which the helicity flux is bounded from below by a positive constant.

An important property of smooth flows of an ideal fluid in two dimensions is conservation of enstrophy (i.e. the $L^2$ norm of the curl of the velocity). In section 4.2 we apply the techniques of Section 3 to the weak formulation of the Euler equations for velocity using a test function that permits estimation of the enstrophy. We obtain the result that, unlike the cases of the energy and the helicity, the locality in the enstrophy cascade is strong only in the ultraviolet range. In the infrared range there are nonlocal effects. Such ultraviolet locality was predicted by Kraichnan \cite{Kraichnan67} and agrees with numerical and experimental evidence. Furthermore, there are arguments in the physical literature that hold that the enstrophy cascade is not local in the infrared range. We present a concrete example that exhibits this behavior.

In the final section of this paper, we study the bilinear term $B(u,v)$. We show that the trilinear map $(u,v,w) \to \langle B(u,v),w) \rangle$ defined for smooth vector fields in $L^3$ has a unique
continuous extension to $\{B^{1/2}_{18/7,2}\}^3$ (and  a fortiori  to $\{H^{5/6}\}^3$, which is the relevant space for the dyadic model problem referred to above). We present an example to show that this result is optimal.
We stress that the borderline space for energy conservation
is much rougher than the space of continuity for $\langle B(u,v), w \rangle$.


\section{Preliminaries}
We will use the notation $\lambda_q = 2^q$ (in some inverse length
units). Let $B(0,r)$ denote the ball centered at $0$ of radius $r$
in $\RR^{d}$. We fix a nonnegative radial function $\chi$ belonging
to~${C_0^{\infty}} (B(0,1))$ such that $\chi(\xi)=1$ for $|\xi|\leq
1/2$. We further define
\begin{equation} \label{defvf}
\vf(\xi) = \chi(\l_1^{-1}\xi) - \chi(\xi).
\end{equation}
Then the following is true
\begin{equation}
\label{lpfond1} \chi(\xi) + \sum_{q\geq 0} \vf (\l_q^{-1}\xi) = 1,
\end{equation}
and
\begin{equation}
\label{lpfond2} |p-q|\geq 2 \Rightarrow \Supp\
\vf(\l_q^{-1}\cdot)\cap \Supp\ \vf(\l_p^{-1}\cdot)=\emptyset.
\end{equation}

We define a Littlewood-Paley decomposition. Let us denote
by $\cF$ the Fourier transform on~$\RR^d$. Let
$h,\
\wt h,\  \D_q $ ($q \ge -1$) be defined as follows:

\begin{align*}
h &= {\cF}^{-1}\vf\quad {\rm and}\quad \wt h =
{\cF}^{-1}\chi, \\
\D_qu & = \cF^{-1}(\vf(\l_q^{-1}\xi)\cF u) = \l_q^d\int h(\l_q y)u(x-y)dy, \, q\ge 0\\
\D_{-1}u & =\cF^{-1}(\chi(\xi)\cF u) = \int \wt h(y)u(x-y)dy.
\end{align*}
For $\L \in \NN$ we define
\begin{equation}
S_\L = \sum_{q=-1}^\L \D_q.
\end{equation}
Due to \eqref{defvf} we have
\begin{equation}
S_\L u = \cF^{-1}(\chi(\l_{\L+1}^{-1} \xi)\cF u).
\end{equation}
Let us now recall the definition of inhomogeneous Besov spaces.

\begin{defi}
Let $s$ be a real number, p and r two real numbers greater than~$1$.
Then
$$
\|u\|_{B^s_{p,r}} \eqdefa \|\D_{-1}u\|_{L^p}+\Bigl\|\left( \l_q^s \|\Delta
_qu\|_{L^p}\right)_{q\in \NN}\Bigr\|_{\ell^r(\NN)}
$$
is the inhomogeneous Besov norm.
\end{defi}

\begin{defi}
Let $s$ be a real number, p and r two real numbers greater than~$1$.
The inhomogeneous Besov space ~$B^s_{p,r}$ is the space of tempered distributions~$u$ such
that the norm ~$\|u\|_{B^s_{p,r}}$ is finite.
\item
\end{defi}
We refer to~\cite{Chemin99} and  \cite{lemarie} for background on
harmonic analysis in the context of fluids. We will use the Bernstein
inequalities
\begin{lemma}
\label{loc}
    $$  \| \D_q u \|_{L^b} \leq \l_q^{d({\frac{1}{ a}} -{\frac{1}{b}})} \|\D_q u
\|_{L^a}
\quad{\rm for }\ b\geq a \geq 1.
$$
\end{lemma}
As a consiquence we have the following inclusions.
\begin{cor}
\label{inclusiontypesob}
{\sl
If~$b\geq a\geq 1$, then we have the following continuous embeddings
\[
B^s_{a,r}\subset B^{s-d\Bigl( \frac 1 a-\frac 1 b\Bigr)}_{b,r}
\]
}
$$
B^0_{a, 2} \subset L^a
$$
\end{cor}
In particular, the following chain of inclusions will be used throughout the text.
\begin{equation}
H^{\frac{5}{6}}(\RR^3) \subset B^{\frac{2}{3}}_{\frac{9}{4},2}(\RR^3) \subset B^{\frac{1}{2}}_{\frac{18}{7},2}(\RR^3)\subset
B^{\frac{1}{3}}_{3,2}(\RR^3).
\end{equation}

\section{Energy flux and locality}

\subsection{Weak solutions}
\begin{defi}
A function $u$ is a weak solution of the Euler equations
with initial data $u_0\in L^2(\RR^d)$
if $u\in C_w([0,T]; L^2(\RR^d))$, (the space of weakly continuous functions)
and for every $\psi \in C^1([0,T]; {\mathcal S}(\RR^d))$ with ${\mathcal S}({\RR^d})$ the space of rapidly decaying functions,
with $\n_x \cdot \psi =0$ and $0\leq t \leq T$, we have
\begin{equation}\label{weaksol}
    (u(t), \psi(t)) - (u(0), \psi(0)) - \int_0^t (u(s),
    \p_s \psi(s) )ds = \int_0^t b(u,\psi,u)(s) ds,
\end{equation}
where
$$
(u,v) = \int_{\RR^d} u\cdot v dx,
$$
$$
b(u,v,w) = \int_{\RR^d} u \cdot \n v \cdot w\ dx,
$$
and $\n_x \cdot u(t) = 0$ in the sense of distributions for every
$t \in [0,T]$.
\end{defi}

Clearly, \eqref{weaksol} implies
Lipschitz continuity of the maps $t \ra (u(t), \psi)$ for fixed test functions. By an approximation argument one can show that for any
weak solution $u$ of the Euler equation, the relationship
\eqref{weaksol} holds for all $\psi$  that are smooth and localized in space, but only weakly Lipschitz in time. This justifies the use of physical space mollifications of $u$ as test functions $\psi$. Because we do not have an existence theory of weak solutions, this is a rather academic point.

\subsection{Energy flux} \label{s:EnergyFlux}

For a divergence-free vector field  $u \in L^2$ we introduce
the Littlewood-Paley energy flux at wave number $\l_\L$ by
\begin{equation}\label{flux}
\Pi_\L = \int_{\RR^3} \tr [ S_\L (u \otimes u) \cdot \n S_\L u] dx.
\end{equation}
If $u(t)$ is a weak solution to the Euler equation, then
substituting the test function $\psi = S^2_\L u$ into the weak formulation
of the Euler equation \eqref{weaksol} we obtain
\begin{equation}\label{solflux}
\Pi_\L(t) = \frac{1}{2} \frac{d }{dt}\|S_\L u(t)\|_2^2.
\end{equation}

Let us introduce the following localization kernel
\begin{equation}
K(q) = \left\{
  \begin{array}{ll}
    \l_q^{2/3}, &  q \leq 0; \\
    \l_q^{-4/3}, & q>0,
  \end{array}
\right.
\end{equation}
For a tempered distribution $u$ in $\RR^3$ we denote
\begin{align}
 d_q &= \l_q^{1/3} \|\D_q u\|_3,\\
 d^2 & = \{d_q^2\}_{q \geq -1}.
\end{align}

\begin{prop}\label{propflux}
The energy flux of a divergence-free vector field  $u \in L^2$
satisfies the following estimate
\begin{equation}\label{fluxloc}
 |\Pi_\L| \leq C (K \ast d^2)^{3/2}(\L).
\end{equation}
\end{prop}
From \eqref{fluxloc} we immediately obtain
\begin{equation}\label{limflux}
 \limsup_{\L \ra \infty} |\Pi_{\L}| \leq \limsup_{\L\ra \infty}d^3_\L.
\end{equation}
We define $\van$ to be  the class of all tempered distributions $u$ in
$\RR^3$ for which
\begin{equation}\label{van}
     \lim_{q \ra \infty} \l_q^{1/3} \|\D_q u\|_3 = 0,
\end{equation}
and hence $d_q \ra 0$. We endow $\van$ with the norm inherited from
$\bes$. Notice that the Besov spaces $B^{1/3}_{3,p}$ for $1\leq p < \infty$,
and in particular $B^{1/3}_{3,2}$ are included in $\van$.

As a consequence of \eqref{solflux} and \eqref{limflux} we obtain the following theorem.

\begin{thm}\label{main}
The total energy flux of any divergence-free vector field in the class
$\van \cap L^2$ vanishes. In particular, every weak solution to the Euler equation
that belongs to the class $L^3([0,T]; \van) \cap C_w([0,T]; L^2)$ conserves energy.
\end{thm}

\begin{proof}[Proof of Proposition \ref{propflux}]
In the argument below all the inequalities should be understood up to
a constant multiple.

Following \cite{cet} we write
\begin{equation}\label{newS}
 S_\L(u \otimes u) = r_\L(u,u) - (u - S_\L)\otimes(u - S_\L) + S_\L u \otimes S_\L u,
\end{equation}
where
\begin{align*}
r_\L(u,u) &= \int_{\RR^3} h_\L(y) (u(x-y) - u(x)) \otimes (u(x-y) - u(x)) dy, \\
\tilde{h}_\L(y)& = \l_{\L+1}^3 \tilde{h}(\l_{\L+1} y).
\end{align*}
After substituting \eqref{newS} into \eqref{flux} we find
\begin{align}
 \Pi_\L & = \int_{\RR^3} \tr [ r_\L(u,u) \cdot \n S_\L u] dx \label{intr}\\
& - \int_{\RR^3} \tr [ (u - S_\L)\otimes(u - S_\L) \cdot \n S_\L u] dx.\label{intS}
\end{align}
We can estimate the term in \eqref{intr} using the H\"{o}lder inequality by
$$
\|r_\L(u,u)\|_{3/2} \|\n S_\L u\|_3,
$$
whereas
\begin{equation}
 \|r_\L(u,u)\|_{3/2} \leq \int_{\RR^3} \left |\tilde{h}_\L(y) \right |\|u(\cdot - y) - u(\cdot)\|_3^2 dy.
\end{equation}
Let us now use Bernstein's inequalities and Corollary \ref{inclusiontypesob} to estimate
\begin{align}
 \|u(\cdot - y) - u(\cdot)\|_3^2 & \leq \sum_{q \leq \L} |y|^2 \l_q^2 \|\D_q u\|_3^2 + \sum_{q>\L} \|\D_q u\|_3^2 \\
& = \l_\L^{4/3} |y|^2 \sum_{q \leq \L} \l_{\L-q}^{-4/3} d_q^2 + \l_{\L}^{-2/3} \sum_{q > \L} \l_{\L - q}^{2/3}
d_q^2 \\
&\leq (\l_\L^{4/3} |y|^2 + \l_{\L}^{-2/3}) (K \ast d^2)(\L).
\end{align}
Collecting the obtained estimates we find
\begin{align*}
 &\left|\int_{\RR^3} \tr [ r_\L(u,u) \cdot \n S_\L u] dx\right|
\\ &\leq (K \ast d^2)(\L)
 \left(\int_{\RR^3} \left |\tilde{h}_\L(y)\right |\l_\L^{4/3} |y|^2 dy + \l_{\L}^{-2/3} \right) \left[\sum_{q \leq \L} \l_q^2 \|\D_q u\|_3^2 \right]^{1/2} \\
 &\leq (K \ast d^2)(\L) \l_{\L}^{-2/3} \left[ \sum_{q \leq \L} \l_q^{4/3} d_q^2 \right]^{1/2} \\
 &\leq (K \ast d^2)^{3/2}(\L)
\end{align*}

Analogously we estimate the term in \eqref{intS}
\begin{align*}
 &\int_{\RR^3} \tr [ (u - S_\L)\otimes(u - S_\L) \cdot \D S_\L u] dx \\& \leq \|u - S_\L u\|_3^2 \|\D S_\L u\|_3\\
& \leq \left( \sum_{q > \L} \|\D_q u\|_3^2 \right) \left( \sum_{q \leq \L} \l_q^2 \|\D_q u\|_3^2 \right)^{1/2}\\
& \leq (K \ast d^2 )^{3/2}(\L).
\end{align*}

This finishes the proof.

\end{proof}

\subsection{Energy flux through dyadic shells}
Let us introduce the energy flux through a sequence of dyadic shells between
scales $-1 \leq \L_0 < \L_1 <\infty$ as follows
\begin{equation}\label{fluxshell}
\Pi_{\L_0\L_1} = \int_{\RR^3} \tr [ S_{\L_0\L_1} (u \otimes u) \cdot \n S_{\L_0\L_1} u] \, dx,
\end{equation}
where
\begin{equation}
 S_{\L_0\L_1} = \sum_{\L_0 \leq q \leq \L_1} \D_q = S_{\L_1} - S_{\L_0}.
\end{equation}
We will show that similar to formula \eqref{fluxloc} the flux through dyadic shells
is essentially controlled by scales near the inner and outer radii. In fact it almost follows from
\eqref{fluxloc} in view of the following decomposition
\begin{equation}
\begin{split}
S_{\L_0\L_1}^2 &= (S_{\L_1} - S_{\L_0-1})^2\\
&=S_{\L_1}^2 + S_{\L_0-1}^2 - 2S_{\L_0-1}S_{\L_1}\\
&=S_{\L_1}^2 + S_{\L_0-1}^2 - 2S_{\L_0-1}\\
&=S_{\L_1}^2 - S_{\L_0-1}^2 -2S_{\L_0-1}(1- S_{\L_0-1})\\
&=S_{\L_1}^2 - S_{\L_0-1}^2 -2\D_{\L_0-1}\D_{\L_0}.
\end{split}
\end{equation}
Therefore
\begin{equation}\label{dsh3}
\Pi_{\L_0\L_1} = \Pi_{\L_1} - \Pi_{\L_0 - 1}
 - 2\int_{\RR^3} \tr [\bar{\D}_{\L_0}(u \otimes u) \cdot \n \bar{\D}_{\L_0} u] \, dx,
\end{equation}
where
\begin{equation}
\bar{\D}_{\L_0}(u) =  \int_{\RR^3} \bar{h}_{\L_0}(y)u(x-y) \, dy,
\end{equation}
and $\bar{h}_{\L_0}(x)= \cF^{-1}\sqrt{\vf(\l^{-1}_{\L_0-1}\xi)\vf(\l^{-1}_{\L_0}\xi)}$.

Note that the flux through a  sequence of dyadic shells is equal to the difference between the
fluxes across the dyadic spheres on the boundary plus a small error term that can be easily
estimated. Indeed, let us rewrite the tensor product term as follows
\begin{equation}
\bar{\D}_{\L_0} (u \otimes u) = \bar{r}_{\L_0}(u,u) + \bar{\D}_{\L_0} u \otimes u+  u \otimes \bar{\D}_{\L_0} u,
\end{equation}
where
$$
\bar{r}_\L(u,u) = \int_{\RR^3} \bar{h}_\L(y) (u(x-y) - u(x)) \otimes (u(x-y) - u(x)) \, dy.
$$
Thus we have
\begin{multline*}
\int_{\RR^3} \tr [ \bar{\D}_{\L_0} (u \otimes u) \cdot \n \bar{\D}_{\L_0} u] \, dx = \int_{\RR^3} \tr [
\bar{r}_\L(u,u) \cdot \n \bar{\D}_{\L_0} u] \, dx \\
- \int_{\RR^3} \bar{\D}_{\L_0} u \cdot \n u \cdot  \bar{\D}_{\L_0} u \, dx
\end{multline*}
We estimate the first integral as previously to obtain
\begin{equation}
 \left|\int_{\RR^3} \tr [ \bar{r}_{\L_0}(u,u) \cdot \n \bar{\D}_{\L_0} u] \, dx \right| \leq
d_{\L_0} (K \ast d^2)(\L_0).
\end{equation}
As to the second integral we have
\begin{equation}
\begin{split}
\left| \int_{\RR^3} \bar{\D}_{\L_0} u \cdot \n u \cdot  \bar{\D}_{\L_0} u \, dx\right| &=
\left| \int_{\RR^3} \bar{\D}_{\L_0} u \cdot \n S_{\L_0} u \cdot  \bar{\D}_{\L_0} u \, dx \right|\\
&\leq d_{\L_0}^2 (K \ast d^2)^{1/2}(\L_0).
\end{split}
\end{equation}

Applying these estimates to the flux \eqref{dsh3} we arrive at the following conclusion.

\begin{thm}
 The energy flux through dyadic shells between wavenumbers $\l_{Q_0}$ and $\l_{Q_1}$ is controlled
primarily by the end-point scales. More precisely, the following estimate holds
\begin{equation}
 |\Pi_{Q_0 Q_1}| \leq C(K\ast d^2)^{3/2}(Q_0) + C(K\ast d^2)^{3/2}(Q_1).
\end{equation}
\end{thm}

\subsection{Construction of a divergence free vector field with
non-vanishing energy flux} \label{energyfluxexample}
In this section we give a construction of a divergence free vector field
in $B^{1/3}_{3,\infty}(\RR^3)$ for which the energy flux
is bounded from below by a positive constant. This suggests
the sharpness of $B^{1/3}_{3,c(\NN)}(\RR^3)$ for energy conservation.
Our construction is based on Eyink's example on a torus \cite{E}, which
we transform to $\RR^3$ using a method described below.

Let $\chi_\L(\xi)=\chi(\lambda_{\L+1}^{-1} \xi)$. We define $P^{\perp}_{\xi}$ for vectors $\xi\in {\RR}^3,\,\,\xi \neq 0$ by
$$
P^{\perp}_{\xi} v = v- |\xi|^{-2}(v\cdot \xi)\xi = \left ({\mathbb I} - |\xi|^{-2}(\xi\otimes\xi)\right)v
$$
for $v\in {\mathbb C}^3$ and we use $v\cdot w =\Sum_{j=1}^3 v_jw_j$
for $v,w\in {\mathbb C}^3$.

\begin{lemma} \label{l:local}
Let $\Phi_k(x)$ be $\RR^3$ -- valued functions, such that
\[I_k := \int_{\RR^3} |\xi||\cF \Phi_k(\xi)| \, d\xi < \infty.\] Let also
$\Psi_k(x)=\mathbb{P}(e^{ik\cdot x} \Phi_k(x))$ where $\mathbb{P}$ is the Leray projector onto the space of divergence free vectors.
Then
\begin{equation}
\sup_x\left|\Psi_k(x) -e^{ik\cdot x}(P^{\perp}_{k} \Phi_{k})(x) \right| \le
\frac{1}{4\pi^3}\frac{I_k}{|k|},
\end{equation}
and
\begin{equation}
\sup_{x}\left |(S^2_\L \Psi_k)(x) - \chi^2_\L(k) \Psi_k(x)\right | \le  \frac{c}{(2\pi)^3}
\frac{I_k}{\lambda_{\L+1}},
\end{equation}
where $c$ is the the Lipschitz constant of $\chi(\xi)^2$.
\end{lemma}
\begin{proof}
First, note that for any $k, \xi \in \RR^3$ and $v \in \CC^3$ we have
\begin{equation} \label{trig1}
\begin{split}
\left| \frac{(v \cdot \xi)\xi}{|\xi|^2} + \frac{(v \cdot \xi)k}{|k|^2} \right|
&\leq \frac{|v|}{|k|}\left| \frac{|k|}{|\xi|}\xi + \frac{|\xi|}{|k|}k \right|\\
&= \frac{|v||\xi+k|}{|k|}.
\end{split}
\end{equation}
In addition, it follows that
\begin{equation} \label{trig2}
\begin{split}
\left| \frac{(v \cdot k)k}{|k|^2} + \frac{(v \cdot \xi)k}{|k|^2} \right|
&= \frac{|(v \cdot (k+\xi))k|}{|k|^2} \\
&\leq \frac{|v||\xi+k|}{|k|}.
\end{split}
\end{equation}
Adding \eqref{trig1} and \eqref{trig2} we obtain
\begin{equation}
\begin{split}
|P^\perp_\xi v - P^\perp_k v| & = \left| \frac{(v \cdot \xi)\xi}{|\xi|^2} - \frac{(v \cdot k)k}{|k|^2} \right|\\
&\leq \left| \frac{(v \cdot \xi)\xi}{|\xi|^2} + \frac{(v \cdot \xi)k}{|k|^2} \right|
+\left| \frac{(v \cdot k)k}{|k|^2} + \frac{(v \cdot \xi)k}{|k|^2}  \right|\\
&\leq 2\frac{|v||\xi+k|}{|k|}.
\end{split}
\end{equation}
Using this inequality we can now derive the following estimate:
\begin{equation}
\begin{split}
|\Psi_k(x) - e^{ik\cdot x}(P^{\perp}_{k} \Phi_k)(x)|&=
|\cF^{-1}[P^{\perp}_{\xi}(\cF \Phi_k)(\xi+k)-P^{\perp}_k(\cF \Phi_k)(\xi+k)]|\\
&\leq  \frac{1}{(2\pi)^3}\int_{\RR^3} 2\frac{|\xi+k|}{|k|} |(\cF \Phi_k)(\xi+k)| \, d\xi \\
&= |k|^{-1}\frac{1}{4\pi^3} \int_{\RR^3} |\xi||(\cF \Phi_k(\xi))| \, d\xi.
\end{split}
\end{equation}
Finally, we have
\begin{equation}
\begin{split}
|(S^2_\L \Psi_k)(x) - \chi_\L(k)^2\Psi_k(x)| &=
|\cF^{-1}[(\chi_\L(\xi)^2-\chi_\L(k)^2)(\cF\Psi_k)(\xi)]|\\
&\leq \frac{1}{(2\pi)^3}\int_{\RR^3} \frac{c|\xi+k|}{\lambda_{\L+1}} |(\cF \Phi_k)(\xi+k)| \, d\xi \\
&= \lambda_{\L+1}^{-1}\frac{c}{(2\pi)^3}\int_{\RR^3} |\xi| |(\cF \Phi_k)(\xi)| \, d\xi,
\end{split}
\end{equation}
where $c$ is the the Lipschitz constant of $\chi(\xi)^2$. This concludes the proof.
\end{proof}

\subsubsection*{Example illustrating the sharpness of Theorem~\ref{main}}

Now we proceed to construct a divergence free vector field
in $B^{1/3}_{3,\infty}(\RR^3)$ with non-vanishing energy flux.
Let $\vu(\vk)$ be a vector field $\vu:\mathbb{Z}^3 \to \mathbb{C}^3$  as in Eyink's example \cite{E}
with
\begin{align*}
U(\lambda_q,0,0)&=i\lambda_q^{-1/3}(0,0,-1),& U(-\lambda_q,0,0)&=i\lambda_q^{-1/3}(0,0,1),\\
U(0,\lambda_q,0)&=i\lambda_q^{-1/3}(1,0,1),& U(0,-\lambda_q,0)&=i\lambda_q^{-1/3}(-1,0,-1),\\
U(\lambda_q,\lambda_q,0)&=i\lambda_q^{-1/3}(0,0,1),& U(-\lambda_q,-\lambda_q,0)&=i\lambda_q^{-1/3}(0,0,-1),\\
U(\lambda_q,-\lambda_q,0)&=i\lambda_q^{-1/3}(1,1,-1),& U(-\lambda_q,\lambda_q,0)&=i\lambda_q^{-1/3}(-1,-1,1),
\end{align*}
for all $q \in \mathbb{N}$ and zero otherwise.
Denote $\rho(x) = \cF^{-1}\chi(4\xi)$ and $A=\int_{\mathbb{R}^3} \rho(x)^3 \, dx$.
Since $\chi(\xi)$ is radial, $\rho(x)$ is real. Moreover,
\begin{align*}
A &= \int_{\mathbb{R}^3} \rho(x)^3 \, dx = \frac{1}{(2\pi)^3} \int_{\mathbb{R}^3} \cF(\rho^2)\cF \rho \, d\xi\\
&= \frac{1}{(2\pi)^6} \int_{\mathbb{R}^3}\int_{\mathbb{R}^3} \chi(4\eta) \chi(4(\xi-\eta)) \chi(4\xi) \, d\eta d\xi >0.
\end{align*}
Now let
\[
u(x) = \mathbb{P}\sum_{k \in \ZZ^3} U(k)e^{ik\cdot x} \rho(x).
\]
Note that $u \in B^{1/3}_{3,\infty}(\RR^3)$. Our goal is to estimate the
flux $\Pi_\L$ for the vector field $u$. Define
\[
\Phi_k = |k|^{1/3} U(k) \rho(x) \qquad \text{and} \qquad \Psi_k(x)=\mathbb{P}(e^{ik\cdot x}\Phi_k(x)).
\]
Then clearly $\Phi_k(x)$ and $\Psi_k(x)$ satisfy the conditions of Lemma~\ref{l:local},
and we have
\begin{equation}
u(x) = \sum_{k\in \ZZ^3 \setminus \{0\}} |k|^{-1/3} \Psi_k(x).
\end{equation}
Now note that
\begin{equation}
\begin{split}
\Psi_{k_1} \cdot \nabla S^2_\L \Psi_{k_2} &=
\Psi_{k_1} \cdot S^2_\L \mathbb{P}[\nabla(e^{ik\cdot x}\Phi_{k_2})]\\
& = i (\Psi_{k_1} \cdot k_2) S^2_\L \Psi_{k_2} + \Psi_{k_1}\cdot S^2_\L \mathbb{P}(e^{ik_2 \cdot x}
\nabla \Phi_{k_2}).
\end{split}
\end{equation}
In addition, the following equality holds by construction:
\begin{equation}
P^{\perp}_{k} \Phi_{k} =  \Phi_{k}, \qquad \forall k \in \ZZ^3.
\end{equation}
Define the annulus $A_\L= \ZZ^3 \cap B(0,\lambda_{\L+2}) \setminus
B(0,\lambda_{\L-1})$.
Thanks to Lemma~\ref{l:local}, for any sequences $k_1(\L), k_2(\L), k_3(\L) \in A_\L$ with $k_1+k_2=k_3$, we have
\[
\begin{split}
\int_{\RR^3} (\Psi_{k_1} \cdot \nabla S^2_\L \Psi_{k_2})\cdot \Psi_{k_3}^* \,dx
&= i\int_{\RR^3} (\Psi_{k_1} \cdot k_2) S^2_\L \Psi_{k_2} \cdot \Psi_{k_3}^* \, dx
+ O(\lambda_{\L}^{0})\\
&\hspace{-1in}=i\int_{\RR^3} (e^{ik_1\cdot x} \Phi_{k_1} \cdot k_2)
\chi_\L(k_2)^2 e^{ik_2\cdot x} \Phi_{k_2} \cdot e^{-ik_3\cdot x} \Phi_{k_3}^* \, dx
+ O(\lambda_{\L}^{0})\\
&\hspace{-1in}=
i(|k_1||k_2||k_3|)^{1/3}A(U(k_1)\cdot k_2)\chi_\L(k_2)^2U(k_2)\cdot U(k_3)^* + O(\lambda_{\L}^{0}).
\end{split}
\]
On the other hand, since the Fourier transform of $\Psi_k$ is supported in \\ $B(k,1/4)$, we have
\begin{equation} \label{e:nointer}
\int_{\RR^3} (\Psi_{k_1} \cdot \nabla S^2_\L \Psi_{k_2})\cdot \Psi_{k_3}^* \,dx
=0,
\end{equation}
whenever $k_1+k_2 \ne k_3$. In addition, due to locality of interactions in this
example, \eqref{e:nointer} also holds if
$A_q \setminus \{k_1,k_2,k_3\} \ne \emptyset$ for all $q \in \NN$. Finally,
\begin{equation}
\int_{\RR^3} (\Psi_{k_1} \cdot \nabla S^2_\L \Psi_{k_2})\cdot \Psi_{k_3}^* \,dx
+\int_{\RR^3} (\Psi_{k_1} \cdot \nabla S^2_\L \Psi_{k_3})\cdot \Psi_{k_2}^* \,dx =0,
\end{equation}
whenever $k_2 \notin A_Q$ and $k_3 \notin A_Q$.
Hence, the flux for $u$ can be written as
\begin{equation}
\Pi_\L =-\sum_{\substack{k_1,k_2, k_3 \in A_\L\\k_1+k_2+k_3=0}}
(|k_1||k_2||k_3|)^{-1/3}
\int_{\RR^3} (\Psi_{k_1} \cdot \nabla S^2_\L \Psi_{k_2})\cdot \Psi_{k_3} \,dx.
\end{equation}
Since the number of nonzero terms in the above sum is independent
of $\L$, we obtain
\begin{equation} \label{errorterm}
\Pi_\L  = A \tilde{\Pi}_\L +O(\lambda_\L^{-1}),
\end{equation}
where $\tilde{\Pi}$ is the flux for the vector field $U$, i.e.,
\begin{equation}
\tilde{\Pi}_{\L} :=
-\sum_{\substack{k_1,k_2, k_3 \in A_\L\\k_1+k_2+k_3=0}}
i(U(k_1)\cdot k_2)\chi_\L(k_2)^2U(k_2)\cdot U(k_3).
\end{equation}
The flux $\tilde{\Pi}_\L$ has only the following
non-zero terms (see \cite{E} for details):
\begin{multline*}
-\sum_{\substack{|k_2|=\lambda_{\L}\\ |k_3|=\sqrt{2}\lambda_\L}} i(\vu_1(-\vk_2-\vk_3) \cdot \vk_2)
\vu_2(\vk_2) \cdot \vu_3(\vk_3) (\chi_\L(k_2)^2-\chi_\L(k_3)^2) \\
\geq 4(\chi(1/2)^2-\chi(1/\sqrt{2})^2),
\end{multline*}
and
\begin{multline*}
-\sum_{\substack{|k_2|=\sqrt{2}\lambda_\L\\|k_3|=2\lambda_\L}} i(\vu_1(-\vk_2-\vk_3) \cdot \vk_2)
\vu_2(\vk_2) \cdot \vu_3(\vk_3)(\chi_\L(k_2)^2-\chi_\L(k_3)^2)\\
\geq 4(\chi(1/\sqrt{2})^2-\chi(1)^2).
\end{multline*}
Hence
\begin{equation*}
\tilde{\Pi}_\L \geq 4(\chi(1/2)^2-\chi(1/\sqrt{2})^2+\chi(1/\sqrt{2})^2-\chi(1)^2) =4.
\end{equation*}
This together with \eqref{errorterm} implies that
\begin{equation*}
\liminf_{\L \to \infty} \Pi_\L \geq 4A.
\end{equation*}

\section{Other conservation laws}

In this section we apply similar techniques to derive optimal results
concerning the conservation of helicity in 3D and that of
enstrophy in 2D for weak solutions of the Euler equation. In the case of the
 helicity flux we prove that simultaneous infrared and ultraviolet
localization occurs, as for the energy flux. However,
 the enstrophy flux exhibits strong
localization only in the ultraviolet region, and a partial
localization in the infrared region. A possibility of such a type of
localization was discussed in \cite{Kraichnan67}.

\subsection{Helicity}
For a divergence-free vector field $u \in H^{1/2}$ with vorticity
$\o = \n \times u \in H^{-1/2}$ we define the helicity and truncated
helicity flux as follows
\begin{align}
\cH &= \int_{\RR^3} u \cdot \o\ dx \\
\cH_\L & = \int_{\RR^3} \tr\left[ S_\L(u \otimes u) \cdot \n S_\L \o
+ S_\L(u \wedge \o) \cdot \n S_\L u  \right] \ dx,
\end{align}
where $u \wedge \o = u \otimes \o - \o \otimes u$. Thus, if $u$ was
a solution to the Euler equation, then $\cH_\L$ would be the time
derivative of the Littlewood-Paley helicity at frequency $\l_\L$,
$$
\int_{\RR^3} S_\L u \cdot S_\L \o\ dx.
$$

Let us denote
\begin{align}
b_q & =\l_q^{2/3} \|\Delta_q u\|_3,\\
b^2 & = \{b_q^2\}_{q=-1}^\infty, \\
T(q) & = \left\{
  \begin{array}{ll}
    \l_q^{2/3}, &  q \leq 0; \\
    \l_q^{-4/3}, & q>0,
  \end{array}
\right.
\end{align}

\begin{prop}\label{helflux}
The helicity flux of a divergence-free vector field  $u \in H^{1/2}$
satisfies the following estimate
\begin{equation}\label{fluxlochel}
 |\cH_\L| \leq C (T \ast b^2)^{3/2}(\L).
\end{equation}
\end{prop}

\begin{thm}\label{mainhel}
The total helicity flux of any divergence-free vector field in the
class $\vanhel \cap H^{1/2}$ vanishes, i.e.
\begin{equation}\label{heltotal}
    \lim_{\L \ra \infty} \cH_\L = 0.
\end{equation}
Consequently, every weak solution to the Euler equation that belongs
to the class $L^3([0,T]; \vanhel) \cap L^\infty([0,T];H^{1/2})$
conserves helicity.
\end{thm}

Proposition~\ref{helflux} and Theorem~\ref{mainhel} are proved by direct
analogy with the proofs of Proposition~\ref{propflux} and Theorem~\ref{main}.

\subsubsection*{Example illustrating the sharpness of Theorem~\ref{mainhel}}

We can also construct an example of a vector field in $B^{2/3}_\infty(\RR^3)$
for which the helicity flux is bounded from below by a positive constant.
Indeed, let $\vu(\vk)$ be a vector field $\vu:\mathbb{Z}^3 \to \mathbb{C}^3$
with
\begin{align*}
U(\pm\lambda_q,0,0)&=\lambda_q^{-2/3}(0,0,-1),\\
U(0,\pm\lambda_q,0)&=\lambda_q^{-2/3}(1,0,1),\\
U(\pm\lambda_q,\pm\lambda_q,0)&=\lambda_q^{-2/3}(0,0,1),\\
U(\pm\lambda_q,\mp\lambda_q,0)&=\lambda_q^{-2/3}(1,1,-1),
\end{align*}
for all $q \in \mathbb{N}$ and zero otherwise.
Denote $\rho(x) = \cF^{-1}\chi(4\xi)$, $A=\int_{\mathbb{R}^3} \rho(x)^3 \, dx$, and
let
\begin{equation}
u(x) = \mathbb{P}\sum_{k \in \ZZ^3} U(k)e^{ik\cdot x} \rho(x).
\end{equation}
Note that $u \in B^{2/3}_{3,\infty}(\RR^3)$. On the other hand, a computation similar to the one
in Section~\ref{energyfluxexample} yields
\begin{equation}
\liminf_{\L\to \infty} \left |\mathcal{H}_\L \right |\geq 4A.
\end{equation}
\subsection{Enstrophy}

We work with the case of a two dimensional fluid in this section. In
order to obtain an expression for the enstrophy flux one can use the
original weak formulation of the Euler equation for velocities
\eqref{weaksol} with the test function chosen to be
\begin{equation}\label{testvort}
    \psi = \n^\perp S_\L^2 \o.
\end{equation}
Let us denote by $\O_\L$ the expression resulting on the right hand
side of \eqref{weaksol}:
\begin{equation}\label{enstflux}
    \O_\L = \int_{\RR^2} \tr\left[ S_\L(u \otimes u) \cdot \n \n^\perp S_\L \o \right] \
    dx.
\end{equation}
Thus,
\begin{equation}\label{fluxder}
    \frac{d \|S_\L \o\|_2^2}{2dt} = \O_\L.
\end{equation}
As before we write
\begin{align*}\label{enstflux2}
    \O_\L & = \int_{\RR^2} \tr\left[ r_\L(u,u) \cdot \n \n^\perp S_\L \o \right] \
    dx \\
    &+ \int_{\RR^2} \tr\left[ (u- S_\L u)\otimes(u - S_\L u) \cdot \n \n^\perp S_\L \o \right] \
    dx
\end{align*}
Let us denote
\begin{align}
c_q & =  \|\Delta_q \o\|_3,\\
c^2 & = \{c_q^2\}_{q=-1}^\infty, \\
W(q) & = \left\{
  \begin{array}{ll}
    \l_q^{2}, &  q \leq 0; \\
    \l_q^{-4}, & q>0,
  \end{array}
\right.
\end{align}

We have the following estimate (absolute constants are omitted)
\begin{align*}
|\O_\L| &\leq \int_{\RR^2} \left |\tilde{h}_\L(y)\right | ( \|\n S_\L u\|_3^2 |y|^2 + \|(I - S_\L)
u\|_3^2 ) \|\n^2 S_\L \o\|_3 dy \\ &+ \|(I-S_\L) u \|_3^2 \|\n^2
S_\L
\o\|_3 \\
& \leq \left(\l_\L^{-2}\|S_\L \o\|_3^2 +\sum_{q > \L} \l_q^{-2}
c_q^2\right)\left(\sum_{q \leq \L} \l_q^4 c_q^2 \right)^{1/2} \\
& + \left( \sum_{q > \L} \l_q^{-2} c_q^2 \right) \left(\sum_{q \leq
\L} \l_q^4 c_q^2 \right)^{1/2} \\
& \leq \|S_\L \o\|_3^2 \left( \sum_{q \leq \L} \l_{\L-q}^{-4} c_q^2
\right)^{1/2} +\left( \sum_{q > \L} \l_{\L-q}^2 c_q^2 \right)\left(
\sum_{q \leq \L} \l_{\L-q}^{-4} c_q^2 \right)^{1/2} \\
& \leq \|S_\L \o\|_3^2 (W \ast c^2)^{1/2}(\L) +(W \ast
c^2)^{3/2}(\L)
\end{align*}

Thus, we have proved the following proposition.

\begin{prop}
The enstrophy flux of a divergence-free vector field satisfies the
following estimate up to multiplication by an absolute constant
\begin{equation}\label{e:enstflux}
|\O_\L| \leq \|S_\L \o\|_3^2 (W \ast c^2)^{1/2}(\L) +(W \ast
c^2)^{3/2}(\L).
\end{equation}
Consequently, every weak solution to the 2D Euler equation \\
with $\o\in L^3([0,T]; L^3)$ conserves enstrophy.
\end{prop}
Much stronger results concerning conservation of enstrophy are
available for the Euler equations (\cite{E05}, \cite{lmn}) and
for the long time zero-viscosity limit for damped and driven Navier-Stokes
equations (\cite{cr}).

\subsubsection*{Example illustrating infrared nonlocality}

We conclude this section with a construction of a vector field for which
the enstrophy cascade is nonlocal in the infrared range.
\begin{figure} \label{f:enstr}
\begin{center}
\includegraphics[width=4in]{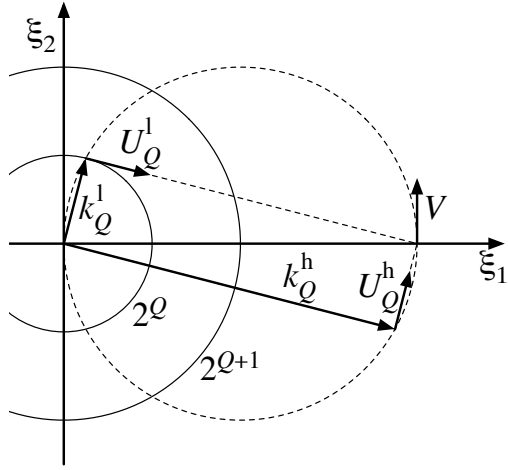}
\end{center}
\caption{Construction of the vector field illustrating infrared nonlocality.}
\end{figure}
Let
$\theta_q = \arcsin(\lambda_{q-\L-2})$ and
\begin{equation}
U^\mathrm{l}_q = (\cos(\theta_q), -\sin(\theta_q)), \qquad
U^\mathrm{h}_q = (\sin(\theta_q), \cos(\theta_q)),
\end{equation}
\begin{equation}
k^\mathrm{l}_q = \lambda_q(\sin(\theta_q), \cos(\theta_q)), \qquad
k^\mathrm{h}_q = \sqrt{\lambda_{\L+2}^2-\lambda_q^2}(\cos(\theta_q), -\sin(\theta_q)),
\end{equation}
see Fig.~\ref{f:enstr} for the case $q=Q$.
Denote $\rho(x) = \delta\tilde{h}(\delta x)$,
$A=\int_{\mathbb{R}^3} \rho(x)^3 \, dx =\int_{\mathbb{R}^3} \tilde{h}(x)^3 \, dx$.
Note that $A>0$ and is independent of $\delta$.
Now let
\begin{equation}
u_q^\mathrm{l}(x) =  \mathbb{P}
[U^\mathrm{l}_q\sin(k^\mathrm{l}_q\cdot x)\rho(x)], \qquad
u_q^\mathrm{h}(x) =  \mathbb{P}[U^\mathrm{h}_q\sin(k^\mathrm{h}_q\cdot x)\rho(x)].
\end{equation}
Let
\begin{equation}
u_q(x)= u_q^\mathrm{l}(x) + u_q^\mathrm{h}(x)
\end{equation}
for $q=0,\dots, \L$, and
\begin{equation}
u_{\L+1}(x)=\mathbb{P}[V\sin(\lambda_{\L+2}x_1)\rho(x)],
\end{equation}
where $V=(0,1)$. Now define
\begin{equation}
u(x) = \sum_{q =0}^{\L+1} u_q(x).
\end{equation}
Our goal is to estimate the enstrophy flux for $u$. Since $\cF u$ is compactly supported,
the expression \eqref{enstflux} is equivalent to
\begin{equation}
\Omega_\L =\int_{\RR^3}(u\cdot \nabla)S^2_\L \omega \cdot \omega \, dx.
\end{equation}
It is easy to see that
\begin{equation}
\Omega_\L \geq \sum_{q=0}^{\L}\int_{\RR^3}(u_q^\mathrm{h}\cdot \nabla)S^2_\L
(\nabla^{\perp} \cdot u_q^{\mathrm{l}}) (\nabla^{\perp} \cdot u_{\L+1}) \, dx.
\end{equation}
Using Lemma~\ref{l:local} we obtain
\begin{equation}
\begin{split}
\Omega_\L  &\geq A\sum_{q=0}^\L |U^\mathrm{h}_q| \lambda_q^2 |U^\mathrm{l}_q| \lambda_{\L+2} |V| + O(\delta)\\
&= \l_{\L+2}\|\D_{\L+2}u\|_3 \sum_{q=0}^\L \l_q^2 \|\D_q u\|_3^2 + O(\delta),
\end{split}
\end{equation}
which shows sharpness of \eqref{e:enstflux} in the infrared range.

\section{Inequalities for the nonlinear term}
We take $d=3$ and consider $u,v\in B^{\frac{1}{3}}_{3,2}$ with $\nabla\cdot u = 0$ and wish to examine the advective term
\be
B(u,v) = {\mathbb P}(u\cdot\nabla v) = \Lambda {\mathbb H}(u\otimes v)
\la{bu}
\ee
where
\be
\left [{\mathbb H}(u\otimes v)\right]_i  = R_j(u_jv_i) + R_i(R_kR_l(u_kv_l))
\la{ha}
\ee
and ${\mathbb P}$ is the Leray-Hodge projector, $\Lambda = (-\D)^{\frac{1}{2}}$ is the Zygmund operator and $R_k = \partial_k\Lambda^{-1}$ are Riesz transforms.

\begin{prop} The bilinear advective term $B(u,v)$ maps continuously
the space $B^{\frac{1}{3}}_{3,2}\times B^{\frac{1}{3}}_{3,2}$ to
the space $B^{-\frac{1}{3}}_{\frac{3}{2}, 2} +
B^{-\frac{2}{3}}_{\frac{9}{5}, 2}$. More precisely, there exist
bilinear continuous maps $C(u,v)$, $I(u,v)$ so that $B(u,v) =
C(u,v) + I(u,v)$ and constants  $C$ such that, for all $u, v \in
B^{\frac{1}{3}}_{3,2}$ with $\nabla\cdot u = 0$, \be
\|C(u,v)\|_{B^{-\frac{1}{3}}_{\frac{3}{2}, 2}} \le
C\|u\|_{B^{\frac{1}{3}}_{3,2}}\|v\|_{B^{\frac{1}{3}}_{3,2}}
\la{inec} \ee and \be \|I(u,v)\|_{B^{-\frac{2}{3}}_{\frac{9}{5},
2}} \le
C\|u\|_{B^{\frac{1}{3}}_{3,2}}\|v\|_{B^{\frac{1}{3}}_{3,2}}\la{inei}
\ee hold. If $u, v, w\in B^{\frac{1}{2}}_{\frac{18}{7}, 2}$ then
\be \left |\langle B(u,v), w\rangle\right | \le
C\|u\|_{B^{\frac{1}{2}}_{\frac{18}{7},2}}\|v\|_{B^{\frac{1}{2}}_{\frac{18}{7},2}}
\|w\|_{B^{\frac{1}{2}}_{\frac{18}{7},2}}\la{tri} \ee holds. So the
trilinear map $(u,v,w) \mapsto \langle B(u,v), w\rangle$ defined
for smooth vector  fields in $L^3$ has a unique continuous
extension to $\left \{B^{\frac{1}{2}}_{\frac{18}{7}, 2}\right
\}^3$ and a fortiori to $\left \{H^{\frac{5}{6}}\right\}^3$.
\end{prop}
\noindent{\bf Proof.} We use duality.
We take $w$ smooth ($w\in B^{\frac{2}{3}}_{\frac{9}{4}, 2}$) and take the scalar product
$$
\langle B(u,v), w\rangle = \int_{\RR^3} B(u,v) \cdot w dx
$$
We write, in the spirit of the paraproduct of Bony
(\cite{Bony81})
\be
\D_q(B(u,v)) = C_q(u,v) + I_q(u,v)\la{bo}
\ee
with
\be
C_q(u,v) = \Sum_{p\ge q-2,\,\, |p-p'|\le 2}\D_q(\Lambda{\mathbb H}(\D_p u, \D_{p'}v))
\la{cq}
\ee
and
\be
I_q(u,v) = \Sum_{j=-2}^2 \left [\D_q\Lambda {\mathbb H}(S_{q+j-2}u, \D_{q+j}v) + \D_q\Lambda {\mathbb H}(S_{q+j-2}v, \D_{q+j}u) \right ]
\la{iq}
\ee
We estimate the contribution coming from the $C_q(u,v)$:
$$
\ba
\Sum_q \left |\langle C_q(u,v), w\rangle \right| \\ \le
C\Sum_{|q-q'|\le 1}\Sum_{p\ge q-2,\,\, |p-p'|\le 2}\lambda_q \lambda_p^{-\frac{2}{3}}\|\Lambda^{\frac{1}{3}}\D_pu\|_{L^3}\|\Lambda^{\frac{1}{3}}\D_{p'}v\|_{L^3}\|\D_{q'}w\|_{L^3} \\
= C\Sum_{|p-p'|\le 2} \lambda_p^{-\frac{2}{3}}\|\Lambda^{\frac{1}{3}}\D_pu\|_{L^3}\|\Lambda^{\frac{1}{3}}\D_{p'}v\|_{L^3}\Sum_{q\le p+2, |q-q'|\le 1}\lambda_q^{\frac{2}{3}}\|\Lambda^{\frac{1}{3}}\D_{q'} w\|_{L^3} \\
\le C\left (\Sum_{|p-p'|\le 2}\|\Lambda^{\frac{1}{3}}\D_pu\|_{\
L^3}\|\Lambda^{\frac{1}{3}}\D_{p'}v\|_{L^3}\right )\|w\|_{B^{\frac{1}{3}}_{3,2}}\\
\le C\|u\|_{B^{\frac{1}{3}}_{3,2}}\|v\|_{B^{\frac{1}{3}}_{3,2}}\|w\|_{B^{\frac{1}{3}}_{3,2}}.
\ea
$$
This shows that the bilinear map $C(u,v) = \sum_{q\ge -1}C_q(u,v)$ maps continuously $\left\{B^{\frac{1}{3}}_{3,2}\right\}^2$ to $B^{-\frac{1}{3}}_{\frac{3}{2}, 2}$ and
\be
\left |\langle C(u,v), w\rangle \right | \le C \|u\|_{B^{\frac{1}{3}}_{3,2}}\|v\|_{B^{\frac{1}{3}}_{3,2}}\|w\|_{B^{\frac{1}{3}}_{3,2}}\la{cineq}
\ee

The terms $I_q(u,v)$ contribute
$$
\ba
\Sum \left |\langle I_q(u,v), w\rangle\right| \\
\le C \Sum_{|j|\le 2, \,\, |q-q'|\le 1}\lambda_q \|S_{q+j-2}u\|_{L^{\frac{9}{2}}}\|\D_{q+j}v\|_{L^3}\|\D_{q'}w\|_{L^{\frac{9}{4}}}
 \\
+ \Sum_{|j|\le 2, \,\, |q-q'|\le 1}\lambda_q\|S_{q+j-2}v\|_{L^{\frac{9}{2}}}\|\D_{q+j}u\|_{L^3}\|\D_{q'}w\|_{L^{\frac{9}{4}}}
\\
\le
C\|u\|_{B^{\frac{1}{3}}_{3,2}}\Sum_{|j|\le 2, |q-q'|\le 1}\lambda_q^{\frac{1}{3}}\|\D_{q+j}v\|_{L^3}\lambda_q^{\frac{2}{3}}\|\D_{q'}w\|_{L^{\frac{9}{4}}} \\
+ C\|v\|_{B^{\frac{1}{3}}_{3,2}}\Sum_{|j|\le 2, |q-q'|\le 1}\lambda_q^{\frac{1}{3}}\|\D_{q+j}u\|_{L^3}\lambda_q^{\frac{2}{3}}\|\D_{q'}w\|_{L^{\frac{9}{4}}}\\
\le C\|u\|_{B^{\frac{1}{3}}_{3,2}}\|v\|_{B^{\frac{1}{3}}_{3,2}}\|w\|_{B^{\frac{2}{3}}_{\frac{9}{4},2}}
\ea
$$
Here we used the fact that
$$
\sup_{q\ge 0}\|S_q u\|_{L^{\frac{9}{2}}} \le C \|u\|_{B^{\frac{1}{3}}_{3,2}}
$$
This last fact is proved easily:
\begin{align*}
\|S_q(u)\|_{L^{\frac{9}{2}}} &\leq \left\|\left(\Sum_{j\le q}|\D_j u|^2\right )^{\frac{1}{2}}\right\|_{L^{\frac{9}{2}}}\\
&\leq \left \{\Sum_{j\le q}\|\D_j u\|_{L^{\frac{9}{2}}}^2\right\}^{\frac{1}{2}}
\leq  C\|u\|_{B^{\frac{1}{3}}_{3,2}}
\end{align*}
We used Minkowski's inequality in $L^{\frac{9}{4}}$ in the penultimate inequality and Bernstein's inequality in the last. This proves that $I$ maps continuously $B^{\frac{1}{3}}_{3,2}\times B^{\frac{1}{3}}_{3,2}$ to $B^{-\frac{2}{3}}_{\frac{9}{5}, 2}$.

The proof of (\ref{tri}) follows along the same lines. Because of Bernstein's inequalities, the inequality (\ref{cineq}) for the trilinear term $\langle C(u,v), w\rangle$ is stronger than (\ref{tri}). The estimate of $I$ follows:
$$
\ba
\Sum \left |\langle I_q(u,v), w\rangle\right| \\
\le C \Sum_{|j|\le 2, \,\, |q-q'|\le 1}\lambda_q \|S_{q+j-2}u\|_{L^{\frac{9}{2}}}\|\D_{q+j}v\|_{L^{\frac{18}{7}}}\|\D_{q'}w\|_{L^{\frac{18}{7}}}
 \\
+ \Sum_{|j|\le 2, \,\, |q-q'|\le 1}\lambda_q\|S_{q+j-2}v\|_{L^{\frac{9}{2}}}\|\D_{q+j}u\|_{L^{\frac{18}{7}}}\|\D_{q'}w\|_{L^{\frac{18}{7}}}\\
\le
C\|u\|_{B^{\frac{1}{3}}_{3,2}}\Sum_{|j|\le 2, |q-q'|\le 1}\lambda_q^{\frac{1}{2}}\|\D_{q+j}v\|_{L^{\frac{18}{7}}}\lambda_q^{\frac{1}{2}}\|\D_{q'}w\|_{L^{\frac{18}{7}}} \\
+ C\|v\|_{B^{\frac{1}{3}}_{3,2}}\Sum_{|j|\le 2, |q-q'|\le 1}\lambda_q^{\frac{1}{2}}\|\D_{q+j}u\|_{L^{\frac{18}{7}}}\lambda_q^{\frac{1}{2}}\|\D_{q'}w\|_{L^{\frac{18}{7}}}\\
\le C\left [\|u\|_{B^{\frac{1}{3}}_{3,2}}\|v\|_{B^{\frac{1}{2}}_{\frac{18}{7},2}} + \|v\|_{B^{\frac{1}{3}}_{3,2}}\|u\|_{B^{\frac{1}{2}}_{\frac{18}{7},2}}\right ]\|w\|_{B^{\frac{1}{2}}_{\frac{18}{7},2}}
\ea
$$
This concludes the proof. $\quad\Box$

The inequality (\ref{cineq}) is not true for $\langle B(u,v), w\rangle$ and (\ref{tri}) is close to being optimal:

\begin{prop} For any $0\le s\le \frac{1}{2}$, $1<p<\infty$, $2<r\le \infty$ there exist functions $u, v, w \in B^{s}_{p,r}$
and smooth, rapidly decaying functions $u_n$, $v_n$, $w_n$, such that
$\lim_{n\to\infty} u_n = u$, $\lim_{n\to\infty}v_n = v$, $\lim_{n\to\infty}w_n=w$ hold in the norm of
$B^{s}_{p,r}$ and such that
$$
\lim_{n\to\infty}\langle B(u_n,v_n), w_n\rangle  = \infty
$$
\end{prop}
\begin{proof}
We start the construction with a divergence-free, smooth function $u$ such that ${\mathcal F}u\in C_0^{\infty}(B(0, \frac{1}{4}))$ and $\int u_1^3dx>0$.
We select a direction $e = (1,0,0)$ and set $\Phi = (0, u_1, 0)$. Then
\be
A := \int_{{\mathbb R}^3} (u(x)\cdot e)\left |P^{\perp}_e \Phi(x)\right |^2dx>0.
\la{aint}
\ee
%
Next we consider the sequence $a_q = \frac{1}{{\sqrt{q}}}$ so that $(a_q) \in \ell^r({\mathbb N})$ for $r>2$, but not for $r=2$, and the functions
\be
v_n = \Sum_{q=1}^{n}\lambda_q^{-\frac{1}{2}}a_q{\mathbb P}\left [\sin(\lambda_qe\cdot x) \Phi(x)\right]
\la{vn}
\ee
and
\be
w_n = \Sum_{q=1}^n \lambda_q^{-\frac{1}{2}}a_q{\mathbb P}\left [\cos(\lambda_q e\cdot x)\Phi(x)\right ].
\la{wn}
\ee
Clearly, the limits $v= \lim_{n\to\infty}v_n$ and $w =\lim_{n\to\infty}w_n$ exist in norm in every $B^{s}_{p, r}$ with $0\le s\le \frac{1}{2}$, $1<p<\infty$ and
$r>2$.
Manifestly, by construction, $u, v_n $ and $w_n$ are divergence-free, and because their Fourier transforms are in $C_0^{\infty}$, they are rapidly decaying functions. Clearly also
$$
\langle B(u, v_n), w_n\rangle = \int_{{\mathbb R}^3}{\mathbb P}(u\cdot\nabla v_n)w_ndx = \int_{{\mathbb R}^3} (u\cdot\nabla v_n)\cdot w_n dx.
$$
The terms corresponding to each $q$ in
\be
\ba
u\cdot\nabla v_n = \Sum_{q=1}^n (u(x)\cdot e)a_q \lambda_q^{\frac{1}{2}}
{\mathbb P}\left[\cos(\lambda_q e\cdot x)\Phi(x)\right] \\
+ \Sum_{q=1}^n a_q\lambda_{q}^{-\frac{1}{2}}u(x)\cdot{\mathbb P}\left[ \sin(\lambda_qe\cdot x)\nabla \Phi(x)\right]
\ea
\la{uv}
\ee
and in (\ref{wn}) have Fourier transforms supported $B(\lambda_q e, \frac{1}{2})\cup B(-\lambda_q, \frac{1}{2})$ and respectively $B(\lambda_q e, \frac{1}{4})\cup B(-\lambda_q e ,\frac{1}{4})$. These are mutually disjoint sets for distinct $q$ and, consequently, the terms corresponding to different indices $q$ do not contribute to the integral $\int(u\cdot\nabla v_n)\cdot w_ndx$. The terms from the second sum in (\ref{uv}) form a convergent series. Therefore, using Lemma~\ref{l:local}, we obtain
\[
\begin{split}
\int_{{\RR}^3}(u\cdot v_n)\cdot w_n &= \Sum_{q=1}^na_q^2\int_{{\RR}^3} (u(x)\cdot e) \left \{\mathbb P\left [\cos(\lambda_q e\cdot x)\Phi(x)\right ]\right\} ^2dx + O(1)\\
&= \Sum_{q=1}^na_q^2\int_{{\RR}^3} (u(x)\cdot e)
\left|P^{\perp}_{e}\Phi(x)\right|^2
\, dx + O(1)\\
&= \left[\Sum_{q=1}^n a_q^2\right] A+O(1),
\end{split}
\]
which concludes the proof.
\end{proof}

\section*{acknowledgment}
The work of AC was partially supported by NSF PHY grant 0555324,
the work of PC  by NSF DMS grant 0504213, the work of SF by NSF DMS grant 
0503768, and the work of RS by NSF DMS grant 0604050.

\end{document}